\documentclass{amsart}

\usepackage{amsmath,amssymb,amscd}

\newtheorem{thm}{Theorem}
\newtheorem{lem}[thm]{Lemma}

\numberwithin{equation}{section}
\numberwithin{thm}{section}

\begin{document}
\title[Noncommutative $L^p$]{Noncommutative $L^p$ structure encodes exactly Jordan structure}
\author{David Sherman}
\address{Department of Mathematics\\
University of Illinois at Urbana-Champaign\\
1409 W. Green Street\\
Urbana, IL 61801-2975}
\curraddr{Department of Mathematics\\ University of California\\ Santa Barbara, CA 93106}
\email{dsherman@math.ucsb.edu}
\urladdr{http://www.math.ucsb.edu/\symbol{126}dsherman}

\begin{abstract}
We prove that for all $1 \le p \le \infty$, $p \ne 2$, the $L^p$ spaces associated to two von Neumann algebras $\mathcal{M}$, $\mathcal{N}$ are isometrically isomorphic if and only if $\mathcal{M}$ and $\mathcal{N}$ are Jordan *-isomorphic.  This follows from a noncommutative $L^p$ Banach-Stone theorem: a specific decomposition for surjective isometries of noncommutative $L^p$ spaces. 
\end{abstract}

\maketitle

\section{Introduction}
In this paper we prove the following theorem.

\begin{thm} \label{T:main}
Let $\mathcal{M}$ and $\mathcal{N}$ be von Neumann algebras, and $1 \le p \le \infty$, $p \ne 2$.  The following are equivalent:
\begin{enumerate}
\item $\mathcal{M}$ and $\mathcal{N}$ are Jordan *-isomorphic;
\item $L^p(\mathcal{M})$ and $L^p(\mathcal{N})$ are isometrically isomorphic as Banach spaces.
\end{enumerate}
\end{thm}

$L^\infty(\mathcal{M})$ is to be understood as $\mathcal{M}$ itself, so for $p=\infty$ the statement follows from the classic article of Kadison \cite{K} (see Theorems \ref{T:jordan} and \ref{T:kadison} below).  One may view this paper as an $L^p$ version of Kadison's results.

The implication (1) $\to$ (2) is a direct application of modular theory and interpolation, only requiring us to go a little further down well-traveled paths.  The more interesting part is to show that (2) $\to$ (1).  In case the surjective isometry is *-preserving and the algebras are $\sigma$-finite, this was proved by Watanabe \cite{W0}.  When both algebras are semifinite, this follows from a structure theorem for $L^p$ isometries (even non-surjective) due to Yeadon \cite{Y}, \cite{Su}; recently Yeadon's theorem was extended in \cite{JRS} to classify isometries for which only the initial algebra is assumed semifinite.  In common with these papers, our proof relies crucially on the equality condition in the noncommutative Clarkson inequality.  But we do not make use of any of these papers' results, and type considerations play no role in our argument (although abelian summands require a little extra care).  We actually determine the structure of the surjective isometry, as follows.

\begin{thm} \label{T:typical} (Noncommutative $L^p$ Banach-Stone theorem)\\
Let $T:L^p(\mathcal{M}) \to L^p(\mathcal{N})$ be a surjective isometry, where $\mathcal{M}$ and $\mathcal{N}$ are von Neumann algebras and $1 < p < \infty$, $p \ne 2$.  Then there are a surjective Jordan *-isomorphism $J:\mathcal{M} \to \mathcal{N}$ and a unitary $w \in \mathcal{N}$ such that
\begin{equation} \label{E:typical}
T(\varphi^{1/p}) = w(\varphi \circ J^{-1})^{1/p}, \qquad \forall \varphi \in \mathcal{M}_*^+.
\end{equation}
\end{thm}

Here $\varphi^{1/p}$ is the generic positive element of $L^p(\mathcal{M})$; we will explain this notation.  Since any $L^p$ element is a linear combination of four positive ones, \eqref{E:typical} completely determines $T$.  The extensions to $0 < p \leq 1$ of Theorems \ref{T:main} and \ref{T:typical} are true but not proved in this paper - see Remark 2 of Section \ref{S:remarks}, and \cite{S3}.

A version of Theorem \ref{T:typical} was shown by Watanabe \cite{W3} under the assumptions that $T$ is *-preserving and $\mathcal{M}$ has a certain extension property.  Our method here is different: we focus on the subspaces $q_1 L^p(\mathcal{M}) q_2$, where $q_1, q_2$ are projections in $\mathcal{M}$.  These subspaces, called corners, are a sort of ``two-dimensional" analogue of the projection bands in classical $L^p$ spaces.  It turns out that $T$ takes corners to corners, preserving both orthogonality (in the sense defined below) and the semi-inner product.  From this we deduce the existence of an orthoisomorphism between the projection lattices of $\mathcal{M}$ and $\mathcal{N}$.  Extending the orthoisomorphism produces a Jordan *-isomorphism, and an intertwining relation finally implies that $T$ has the form \eqref{E:typical}.

Theorem \ref{T:typical} evidently suggests the larger challenge of classifying all $L^p$ isometries.  While this is still open in general, we mention that the author has recently written an article \cite{S1} which obtains several new results, including a solution which is valid under a mild (perhaps vacuous?) hypothesis on the initial algebra.  Also the paper \cite{JRS} completely determines the structure of 2-isometries between $L^p$ spaces.  Although there is some overlap in the setup of these problems, we believe that the surjective case merits a separate exposition, being of independent interest and admitting a distinct technique and solution.  There is no overlap at all - in fact, an interesting contrast - with the investigation \cite{HRS} into non-isometric embeddings between noncommutative $L^p$ spaces.

\section{Background}

We start with some notation.  The only algebras (denoted $\mathcal{M}$, $\mathcal{N}$) under consideration in this paper are von Neumann algebras.  We will use $\mathcal{Z}$ for ``center of" and $\mathcal{P}$ for ``projections of", so for example $\mathcal{P}(\mathcal{Z}(\mathcal{M}))$ is the set of central projections of $\mathcal{M}$.  With $\varphi \in \mathcal{M}_*,$ $x \in \mathcal{M}$, $x\varphi$ (resp. $\varphi x$) means the functional $\varphi(\cdot \: x)$ (resp. $\varphi(x \: \cdot)$).  We use $s_\ell, s_r$ to mean ``left/right support of", for operators, functionals, or $L^p$ vectors.  Often we simply write $L^p$ to indicate a generic noncommutative $L^p$ space.

A \textit{Jordan homomorphism} between von Neumann algebras is a linear map which preserves the Jordan operator product $(x,y) \mapsto (1/2)(xy + yx)$.  Possible adjectives include normal, *-preserving, injective, surjective... a Jordan homomorphism which is all of these is a \textit{surjective Jordan *-isomorphism}.  (Normality is a consequence \cite[Paragraph 4.5.6]{HS}.)  That being said, all of the Jordan theory that the reader needs for this paper is contained in Kadison's
\begin{thm} \label{T:jordan} \cite[Theorem 10]{K}
A surjective Jordan *-isomorphism between von Neumann algebras is the direct sum of a *-isomorphism and a *-anti\-iso\-morphism.
\end{thm}

Up to multiplication by a unitary, these are all the surjective isometries between von Neumann algebras.

\begin{thm} \label{T:kadison} \cite[Theorem 7]{K}
Let $T$ be a surjective isometry between the von Neumann algebras $\mathcal{M}$ and $\mathcal{N}$.  Then there are a surjective Jordan *-isomorphism $J$ from $\mathcal{M}$ to $\mathcal{N}$ and a unitary $w \in \mathcal{N}$ such that $T(x) = wJ(x)$ for all $x \in \mathcal{M}$.
\end{thm}

Actually Kadison proved both of these theorems for all unital C*-algebras.  Since isometries of abelian unital C*-algebras are described by the Banach-Stone theorem, Theorem \ref{T:kadison} is considered a noncommutative Banach-Stone theorem.  The reader will note its similarity with Theorem \ref{T:typical}.  But Kadison's proof of Theorem \ref{T:kadison}, and others offered later, rely on the geometry (i.e. extreme points, faces) of the unit ball.  It does not seem that they can be adapted to work in the $L^p$ context.

We will assume a basic familiarity with noncommutative $L^p$ spaces.  Still, it seems wise to review briefly the specific constructions and concepts that we need.  We provide selective, but hopefully sufficient, references to the literature.  The reader desiring more overview might consult \cite{PX}, which focuses on Banach space properties and also includes a rich bibliography.

\bigskip

In keeping with the motto ``von Neumann algebras are noncommutative $L^\infty$ spaces", one thinks of von Neumann preduals as noncommutative $L^1$ spaces and can consider how to construct their $L^p$ cousins.  When $\mathcal{M}$ is a semifinite algebra with faithful normal semifinite tracial weight $\tau$, one may simply employ $\tau$ as an integral.  That is, $L^p(\mathcal{M}, \tau)$ is the closure of $\{T \in \mathcal{M} \mid \|T\|_p \triangleq \tau(|T|^p)^{1/p} < \infty \}$ in the norm $\|\cdot\|_p$.  This construction goes back to Segal \cite{Se} and has a pleasing interpretation as a set of (possibly unbounded) operators.  See \cite{N}.

But it does not work for all von Neumann algebras.  The first general construction is due to Haagerup \cite{H}, who saw that $\mathcal{M}_*$ could be identified, as an ordered vector space, with a class of unbounded operators affiliated with the core of $\mathcal{M}$.  Since these are operators, one can take $p$th roots on the positive cone, and the norm can be imported from $\mathcal{M}_*$.  To be specific, $L^p(\mathcal{M})$ is the set of $\tau$-measurable operators affiliated with $(\mathcal{M} \rtimes_\sigma \mathbb{R}, \tau)$ which satisfy $\theta_s(T) = e^{-s/p}T$.  Here $\sigma$ is a modular action, $\tau$ is the canonical trace, and $\theta$ is the dual action.  Notice that the product of an $L^p$ operator and an $L^q$ operator is an $L^r$ operator, where $\frac{1}{p} + \frac{1}{q} = \frac{1}{r}$.  See \cite{T1}.

In this construction, any positive element in $L^p(\mathcal{M})$ is the $p$th root of an operator which corresponds to some $\varphi \in \mathcal{M}_*^+$.  We will refer to this element as $\varphi^{1/p}$.  Notice that $\|\varphi^{1/p}\| = [\varphi(1)]^{1/p}$.  This notation frequently proves expedient and is discussed specifically in \cite{Ya}, \cite[Section V.B.$\alpha$]{C2}, \cite{S2}.

The polar decomposition and $\mathcal{M}-\mathcal{M}$ bimodule structure for $L^1(\mathcal{M})$ agree with those of $\mathcal{M}_*$.  In particular, the partial isometry and support projections are in $\mathcal{M}$, and all support projections are necessarily $\sigma$-finite.  This second statement remains true for $L^p(\mathcal{M})$, but the bimodule structure is less obvious.  See \cite{T1}, \cite{JS}.

Another construction of $L^p(\mathcal{M})$ is by complex interpolation, pioneered by Kosaki \cite{Ko1}.  Assume that $\mathcal{M}$ is $\sigma$-finite, and consider the left embedding of $\mathcal{M}$ in $\mathcal{M}_*$ arising from a fixed faithful state $\varphi \in \mathcal{M}_*$: $x \mapsto x\varphi$.  Then Haagerup's space $L^p(\mathcal{M})$ is isometric to the interpolated Banach space at $1/p$ \cite[Theorem 9.1]{Ko1}.  More precisely, we have
\begin{equation} \label{E:kosaki}
L^p(\mathcal{M}) \simeq [\mathcal{M}, \mathcal{M}_*]_{1/p} = L^p(\mathcal{M}) \varphi^{1/q}, \qquad 1/p + 1/q = 1.
\end{equation}
Here the equality is meant as \textit{sets}, while the isomorphism is an isometric identification of Banach spaces.  Right embeddings of the form $x \mapsto \varphi x$ (and even others) work equally well.

Evaluation at 1 (i.e., $\varphi \mapsto \varphi(1)$) is a distinguished linear functional on $\mathcal{M}_* \simeq L^1(\mathcal{M})$.  It is called the \textit{Haagerup trace}, and denoted Tr, because it implements the duality between $L^p$ and $L^q$ ($1/p + 1/q = 1$) in a trace-like way:
$$<\xi, \eta> = \text{Tr}(\xi \eta) = \text{Tr}(\eta \xi), \qquad \xi \in L^p(\mathcal{M}), \: \eta \in L^q(\mathcal{M}).$$
Under this pairing each of $L^p(\mathcal{M})$ and $L^q(\mathcal{M})$ can be isometrically identified with the dual space of the other, and of course $L^\infty(\mathcal{M}) = \mathcal{M}$ is the dual space of $L^1(\mathcal{M})$ \cite{T1}.

The most important $L^p$ result for this paper is the noncommutative Clarkson inequality, or more accurately the condition characterizing when it is an equality.  Yeadon \cite{Y} showed this for semifinite von Neumann algebras; a few years later Kosaki \cite{Ko2} proved it for arbitrary von Neumann algebras with $2<p<\infty$; and only recently Raynaud and Xu \cite{RX} obtained a general version (relying on Kosaki's work).

\begin{thm} \label{T:clarkson} (Equality condition for noncommutative Clarkson inequality)

For $\xi, \eta \in L^p$, $0 < p < \infty$, $p \ne 2$,
\begin{equation} \label{E:clarkson}
\|\xi + \eta\|^p + \|\xi - \eta\|^p = 2(\| \xi \|^p + \| \eta \|^p) \iff \xi \eta^* = \xi^* \eta = 0.
\end{equation}
\end{thm}

The second condition is equivalent to requiring $s_r(\xi) s_r(\eta) = s_\ell(\xi) s_\ell(\eta) = 0$.  Because of this, we call pairs of $L^p$ vectors satisfying \eqref{E:clarkson} \textbf{orthogonal}.  Since the first condition of \eqref{E:clarkson} is preserved by isometries, orthogonality is preserved too.  (For classical $L^p$ spaces, this says that isometries preserve disjointness of support.  Banach made this observation in the very first investigation of $L^p$ isometries \cite{B}.)  To keep things clear, this is the only usage of the term ``orthogonal" in this paper, except where we refer specifically to orthogonality of projections.  We do not use ``orthogonal" to describe pairs of vectors with semi-inner product zero.  So for a set $S \subset L^p$, the orthocomplement $S^\perp$ means the set of elements orthogonal (in this sense) to every element in $S$.

Some authors use ``disjoint" in place of ``orthogonal".  We reserve this term for another use: two subspaces are called \textbf{disjoint} if their intersection is $\{0\}$.

\section{Corners and semi-inner products}

It will be helpful to introduce some \textit{ad hoc} terminology: a subspace of $L^p$ is a \textbf{corner} if it is of the form $q_1 L^p q_2$ for some projections $q_1,q_2$.  Corners with $q_1 = 1$ (resp. $q_2 = 1$) will be called \textbf{columns} (resp. \textbf{rows}).  Note that a corner has a unique representation in which $q_1, q_2$ have equal central support; by the central support of a corner we mean the central support of the projections in such a representation.  We also refer to either $\mathcal{M}z$ or $L^p(\mathcal{M})z$ as a \textbf{central summand} when $z \in \mathcal{P}(\mathcal{Z}(\mathcal{M}))$.  

\begin{lem} \label{T:perp} ${}$
\begin{enumerate}
\item If $T$ is a surjective isometry between $L^p$ spaces (\,$1 \le p < \infty$, $p \ne 2$) and $S$ is a subset of the domain, then $T(S^\perp) = T(S)^\perp$.
\item The intersection of any collection of corners is a corner.
\item For any set $S \subset L^p$, $S^\perp$ is a corner.
\end{enumerate}
\end{lem}

\begin{proof}
$T$ and $T^{-1}$ preserve orthogonality, proving the first statement.  For the second, let $\{p_\alpha\}, \{q_\alpha\}$ be sets of projections; then
$$\bigcap p_\alpha L^p q_\alpha = (\wedge p_\alpha) L^p (\wedge q_\alpha).$$
The third follows from noting that $\{\xi\}^\perp = (1-s_\ell(\xi)) L^p (1-s_r(\xi))$ and applying the second to the expression
$$S^\perp = \bigcap_{\xi \in S} \{\xi\}^\perp. \qedhere$$
\end{proof}

\bigskip

The other notion we need is that of a \textbf{semi-inner product}, first defined for general Banach spaces by Lumer \cite{L}.  We will specialize our discussion to $L^p$ spaces, $1 < p < \infty$.  A nice development of the relationship between isometries and semi-inner products can be found in \cite[Section 1.4]{FJ}.

For $\eta \in L^p$, define $\varphi_\eta$ to be the unique functional in $(L^p)^*$ with $\|\varphi_\eta\| = \|\eta\|$ and $\varphi_\eta(\eta) = \|\eta\|^2$.  The assignment $\eta \mapsto \varphi_\eta$ is known as a \textit{duality map}; uniqueness of the duality map is expressed by saying that $L^p$ is a \textit{smooth} Banach space.  We have that $\varphi_0 = 0$ and otherwise
\begin{equation} \label{E:lpsip}
\varphi_\eta (\cdot) = \frac{\text{Tr} (\cdot \: |\eta|^{p-1} v^*)}{\|\eta\|^{p-2}},
\end{equation}
where $\eta$ has polar decomposition $v|\eta|$.  The semi-inner product is the function on $L^p \times L^p$ defined by
\begin{equation} \label{E:sip}
[\xi, \eta] \triangleq \varphi_\eta(\xi), \qquad \xi, \eta \in L^p.
\end{equation}
In general the semi-inner product is not additive in the second variable.

We prepare two lemmas for later use.  The first is a small variation of well-known results and surely appears in the literature somewhere.  See \cite{KR} for the historical predecessor or \cite[Lemma 4.2]{JRS} for a similar application.

\begin{lem} \label{T:preserve}
If $T$ is an isometry between $L^p$ spaces (\,$1 < p < \infty$), then $T$ preserves the semi-inner product.
\end{lem}

\begin{proof}
Note that we are not assuming that $T$ is surjective, so that $T^*$ is only contractive.  We first take any $L^p$ vector $\xi$ and calculate
$$T^*(\varphi_{T\xi})(\xi) = \varphi_{T\xi}(T \xi) = [T\xi, T\xi] = \|T \xi\|^2 = \|\varphi_{T\xi}\| \|\xi\| \geq \| T^*(\varphi_{T\xi})\| \|\xi\|,$$
so by smoothness $T^*(\varphi_{T\xi}) = \varphi_\xi.$  Now we apply this to any two $L^p$ vectors $\xi, \eta$:
$$[T\xi, T\eta] = \varphi_{T\eta}(T\xi) = T^*(\varphi_{T\eta})(\xi) = \varphi_\eta(\xi) = [\xi, \eta]. \qedhere$$
\end{proof}

\begin{lem} \label{T:corth}
Let $1 < p < \infty$, and let $p_1 L^p p_2$ and $q_1 L^p q_2$ be corners such that 
\begin{equation} \label{E:sipzero}
[\xi, \eta] = 0, \qquad \forall \xi \in p_1 L^p p_2, \: \forall \eta \in q_1 L^p q_2.
\end{equation}
Then $p_1 q_1$ and $p_2 q_2$ are centrally orthogonal. 
\end{lem}

\begin{proof}
Using \eqref{E:lpsip}, \eqref{E:sipzero} is equivalent to
$$\text{Tr}(p_1 \xi p_2 q_2 \zeta q_1) = 0, \qquad \forall \xi \in L^p, \forall \zeta \in L^q, \qquad \frac{1}{p} + \frac{1}{q} = 1.$$
By duality we may conclude that $q_1 p_1 \xi p_2 q_2 = 0$ for any $\xi \in L^p$.  Since the central supports of $q_1 p_1$ and $(q_1 p_1)^* = p_1 q_1$ are equal, this implies the lemma.
\end{proof}

\section{Proof of Theorems \ref{T:main} and \ref{T:typical}}

Let us start with the implication (1) $\to$ (2) of Theorem \ref{T:main}.  The case $p=\infty$ is automatic; Theorem \ref{T:jordan} shows that a surjective Jordan *-isomorphism is isometric.  The case $p=1$ follows by considering the preadjoint of the (normal) surjective Jordan *-isomorphism.  We now assume $1< p<\infty$, $p \ne 2$, and the existence of a surjective Jordan *-isomorphism $J:\mathcal{M} \to \mathcal{N}$.

By Theorem \ref{T:jordan}, there is a central projection $z \in \mathcal{M}$ such that $xz \mapsto J(x)J(z)$ is a surjective *-isomorphism from $z\mathcal{M}$ to $J(z)\mathcal{N}$, and $x(1-z) \mapsto J(x)J(1-z)$ is a surjective *-antiisomorphism from $(1-z)\mathcal{M}$ to $J(1-z)\mathcal{N}$.  Since $L^p(\mathcal{M})$ is isometric to $L^p(z\mathcal{M}) \oplus_p L^p((1-z)\mathcal{M})$ (and similarly for $\mathcal{N}$), it suffices to show that *-isomorphic or *-antiisomorphic von Neumann algebras have isometric $L^p$ spaces.

At least the *-isomorphic case is known.  In fact the core of a von Neumann algebra, so also its $L^p$ spaces, can be constructed \textit{functorially} (see, for example, \cite[Theorem 3.5]{FT}).  Here we cover the *-antiisomorphic case only; the reader will have no trouble making the necessary changes for a *-isomorphism.  A related discussion is in \cite[Section 3]{W1}, although some statements were later corrected in \cite[Section 3]{W2}.

So let $\alpha: \mathcal{M} \to \mathcal{N}$ be a surjective *-antiisomorphism.  (This does not imply that there exists a surjective *-isomorphism, by a paper of Connes \cite{C1}.)  We want to construct a surjective isometry from $L^p(\mathcal{M})$ to $L^p(\mathcal{N})$.

Temporarily assume that the algebras are $\sigma$-finite, and fix a faithful state $\varphi \in \mathcal{M}_*^+$.  We know that $L^p(\mathcal{M}) \simeq [\mathcal{M}, \mathcal{M}_*]_{1/p}$ and $L^p(\mathcal{N}) \simeq [\mathcal{N}, \mathcal{N}_*]_{1/p}$, where we use the embeddings
$$\mathcal{M} \ni x \overset{\iota_1}{\mapsto} x \varphi \in \mathcal{M}_*, \qquad \mathcal{N} \ni y \overset{\iota_2}{\mapsto} (\varphi \circ \alpha^{-1}) y \in \mathcal{N}_*.$$
Then the following diagram commutes, and the horizontal arrows are isometric linear isomorphisms.
\[
\begin{CD}
\mathcal{M} @>\alpha>> \mathcal{N}\\
@V\iota_1VV            @VV\iota_2V\\
\mathcal{M}_*  @>>(\alpha^{-1})_*>  \mathcal{N}_*
\end{CD}
\]

It follows that the interpolated spaces are isometrically isomorphic, and the $\sigma$-finite case is settled.  One might handle the non-$\sigma$-finite case by interpolating with a faithful (normal semifinite) \textit{weight}.  The first $L^p$ construction along these lines is \cite{T2}; \cite{I} marshals even more technical machinery to recover the analogues of the left and right embeddings above.  We will go in a different direction.  

If we look at the equality \eqref{E:kosaki}, we see that $x\varphi^{1/p} \in L^p(\mathcal{M})$ is being identified in $\mathcal{M}_*$ with $x\varphi$.  This corresponds to $(\varphi \circ \alpha^{-1}) \alpha(x) \in \mathcal{N}_*$, which gives the $L^p$ element $(\varphi \circ \alpha^{-1})^{1/p} \alpha(x)$.  So the isometry is densely defined by
$$x\varphi^{1/p} \mapsto (\varphi \circ \alpha^{-1})^{1/p} \alpha(x), \qquad x \in \mathcal{M}.$$

Actually, this map is independent of the choice of $\varphi$.  We have that
\begin{equation} \label{E:ind}
x \varphi^{1/p} = y \psi^{1/p} \Rightarrow (\varphi \circ \alpha^{-1})^{1/p} \alpha(x) = (\psi \circ \alpha^{-1})^{1/p} \alpha(y),
\end{equation}
using the cocycle identity
\begin{equation} \label{E:anticocyc}
(D(\psi \circ \alpha^{-1}): D(\varphi \circ \alpha^{-1}))_t = \alpha((D\varphi:D\psi)_{-t}).
\end{equation}
Equations \eqref{E:ind} and \eqref{E:anticocyc} are checked explicitly in \cite[Section 6]{S1}, based on \cite[Corollary VIII.1.4 and Theorem VIII.3.3]{Ta}.  But this is not yet enough to conclude that the isometries associated to $\varphi$ and $\psi$ are equal, as the subspace $\mathcal{M} \varphi^{1/p} \cap \mathcal{M} \psi^{1/p}$ may not be dense in $L^p(\mathcal{M})$.  (See \cite{Ko3} for a discussion of non-density when $p=2$.)  However, given faithful states $\varphi, \psi \in \mathcal{M}_*$, we may use functional calculus to define the auxiliary state
$$L^1(\mathcal{M}) \ni \rho = \frac{(\varphi^{2/p} + \psi^{2/p})^{p/2}}{\|(\varphi^{2/p} + \psi^{2/p})^{p/2}\|_1}.$$
From the $\tau$-measurable operator inequality $\varphi^{2/p} \le C \rho^{2/p}$, it follows that $\varphi^{1/p} = x \rho^{1/p}$ for some $x \in \mathcal{M}$.  This means that $\mathcal{M}\varphi^{1/p} \cap \mathcal{M}\rho^{1/p} = \mathcal{M}\varphi^{1/p}$, which is dense in $L^p(\mathcal{M})$.  Then \eqref{E:ind} shows that $\varphi$ and $\rho$ generate the same isometry.  But $\psi$ and $\rho$ generate the same isometry too, so in the end we can identify the isometries from $\varphi$ and $\psi$.  Some of the details of this argument are given in \cite[Section 1]{JS}, \cite{S2}, and also generalized in \cite[Section 6]{S1}.

Let us call this $L^p$ isometry $\alpha_p$.  The independence of $\alpha_p$ from any choice of functional implies that 
\begin{equation} \label{E:alphap}
\alpha_p(\varphi^{1/p}) = (\varphi \circ \alpha^{-1})^{1/p}, \qquad \forall \varphi \in \mathcal{M}_*.
\end{equation}
This can actually be taken as a definition for $\alpha_p$, since every element in $L^p(\mathcal{M})$ is a linear combination of four positive ones.  Notice also that $\alpha_p(x \xi) = \alpha_p(\xi) \alpha(x).$  Equation \eqref{E:alphap} tells us that $\alpha_p$ is positive (thus *-preserving), so we can improve this to
\begin{equation} \label{E:bimod}
\alpha_p(x \xi y) = \alpha(y) \alpha_p(\xi) \alpha(x), \qquad \xi \in L^p(\mathcal{M}),\: x,y \in \mathcal{M}.
\end{equation}

Now for any $\sigma$-finite $q \in \mathcal{P}(\mathcal{M})$, we can construct a surjective isometry from $qL^p(\mathcal{M})q$ to $\alpha(q) L^p(\mathcal{N}) \alpha(q)$ as above.  Every finite set of vectors in $L^p(\mathcal{M})$ belongs to some such $qL^p(\mathcal{M})q$, as the left and right supports of each vector belong to the lattice of $\sigma$-finite projections.  Furthermore, these isometries can be defined by \eqref{E:alphap}, so they agree on common domains.  It follows that \eqref{E:alphap} defines a global $L^p$ isometry in the non-$\sigma$-finite case as well.

This ends the proof of (1) $\to$ (2).  More discussion of $L^p$ isometries constructed by interpolation, involving conditional expectations or more general projections, can be found in \cite[Sections 6 and 7]{S1}.

\bigskip

We now turn to the implication (2) $\to$ (1) of Theorem \ref{T:main}.  When $p=\infty$, this follows from Theorem \ref{T:kadison}.  In case $p=1$, the adjoint of a surjective isometry is again a surjective isometry, and we may appeal to the preceding statement.  The implication for the remaining values of $p$ is an obvious consequence of Theorem \ref{T:typical}, which we prove in the remainder of this section.  
Assume that $T: L^p(\mathcal{M}) \to L^p(\mathcal{N})$ is a surjective isometry of Banach spaces, with $1 < p < \infty$, $p \ne 2$.

\begin{lem} \label{T:corner}
If $z \in \mathcal{P}(\mathcal{Z}(\mathcal{M})),$ then
\begin{equation} \label{E:center}
T(z L^p(\mathcal{M})) = z' L^p(\mathcal{N}) \: \text{ for some } z' \in \mathcal{P}(\mathcal{Z}(\mathcal{N})).
\end{equation}
The map $z \mapsto z'$ induces a surjective *-isomorphism from $\mathcal{Z}(\mathcal{M})$ to $\mathcal{Z}(\mathcal{N})$.  
\end{lem}

\begin{proof}
The corners $zL^p(\mathcal{M})$ and $(1-z)L^p(\mathcal{M})$ are orthocomplements of each other, so by Lemma \ref{T:perp}(1) their images are orthocomplements of each other.  Then Lemma \ref{T:perp}(3) tells us there are $q,r,s,t \in \mathcal{P}(\mathcal{N})$ with 
\begin{equation} \label{E:censum}
T(zL^p(\mathcal{M})) = q L^p(\mathcal{N}) r, \qquad T((1-z)L^p(\mathcal{M})) = sL^p(\mathcal{N})t.
\end{equation}
We may assume that the central supports of $q$ and $r$ are equal, and of $s$ and $t$ are equal.  From \eqref{E:censum} it follows that each vector in $L^p(\mathcal{N})$ can be uniquely written as the sum of two orthogonal vectors, one from each of $q L^p(\mathcal{N}) r$ and $s L^p(\mathcal{N}) t$.  As projections, $q,s$ are orthogonal, and $r,t$ are orthogonal.  The conclusion \eqref{E:center} will follow if we can show that $q$ and $t$ are centrally orthogonal projections, for then the spanning property just mentioned implies that all four projections are central.

If $q$ and $t$ are not centrally orthogonal, we can find $0 \ne \xi \in qL^p(\mathcal{N})t$.  Write the decomposition as $\xi = \xi_1 + \xi_2$, and note that neither of $\xi_1, \xi_2$ can be zero.  Now the left support of a sum of orthogonal vectors is the sum of the left supports, just as it is for operators.  So $s_\ell(\xi) \nleq q$, which is a contradiction.

Since $T^{-1}$ also satisfies \eqref{E:center}, the correspondence $z \leftrightarrow z'$ is bijective.  It is additive on orthogonal elements and so induces a surjective *-isomorphism.
\end{proof}

Lemma \ref{T:corner} is related to \cite[Proposition 6.2]{JS}.  In the sequel we use the apostrophe to indicate the correspondence $z \leftrightarrow z'$ without further mention.

\begin{lem} \label{T:abelian}
Let $a \in \mathcal{P}(\mathcal{Z}(\mathcal{M}))$ be such that $a \mathcal{M}$ is the abelian summand of $\mathcal{M}$.  Then $a'\mathcal{N}$ is the abelian summand of $\mathcal{N}$.
\end{lem}

\begin{proof}
We first argue that $a'\mathcal{N}$ is abelian.  If if is not, let $q$ be a noncentral projection in $a'\mathcal{N}$.  Since $qL^p(\mathcal{N})q = [(1-q)L^p(\mathcal{N})(1-q)]^\perp$, we have by Lemma \ref{T:perp} that $T^{-1}(q L^p(\mathcal{N})q)$ is a corner of $L^p(\mathcal{M})$.  But $T^{-1}(q L^p(\mathcal{N})q)$ is contained in $aL^p(\mathcal{M})$, so being a corner it must be a central summand.  Using Lemma \ref{T:corner} we conclude that $q L^p(\mathcal{N})q = T[T^{-1}(q L^p(\mathcal{N})q)]$ is a central summand, which is a contradiction.

Combined with a symmetric argument for $T^{-1}$, this proves the lemma.
\end{proof}

\begin{lem} \label{T:allcorner} ${}$
\begin{enumerate}
\item The image of any corner under $T$ is again a corner.
\item If $q \in \mathcal{P}(\mathcal{M})$ is strictly between 0 and 1 on all central summands, then
\begin{equation} \label{E:tcol}
T(L^p(\mathcal{M})q) = L^p(\mathcal{N}) q_1 z' + q_2 L^p(\mathcal{N})(1-z'), \end{equation}
for some $q_1,q_2 \in \mathcal{P}(\mathcal{N}),$ $z' \in \mathcal{P}(\mathcal{Z}(\mathcal{N})),$ with $q_1 z' + q_2 (1-z')$ strictly between 0 and 1 on every central summand.
\end{enumerate}
\end{lem}

\begin{proof}
For the first statement, let $p_1 L^p(\mathcal{M}) p_2$ be an arbitrary corner.  Then there are central projections $y_1,y_2,y_3,y_4$ with sum $1$, such that
\begin{itemize}
\item $p_1, p_2$ are strictly between 0 and 1 on every central subsummand of $\mathcal{M}y_1$; 
\item $p_1 L^p(\mathcal{M}) p_2 y_2$ is a column which contains no central summand and has central support $y_2$;
\item $p_1 L^p(\mathcal{M}) p_2 y_3$ is a row which contains no central summand and has central support $y_3$;
\item $p_1 L^p(\mathcal{M}) p_2 y_4$ is a central summand.
\end{itemize}

By Lemma \ref{T:corner}, $T$ preserves central sums and takes central summands to central summands.  Therefore we may treat each of the cases separately, and the fourth case is clear.  For the first case, $p_1 L^p(\mathcal{M}) p_2 y_1$ and $(1- p_1) L^p(\mathcal{M})(1-p_2) y_1$ are orthocomplements in $L^p(\mathcal{M} y_1)$, so by Lemma \ref{T:perp} their images are corners.  The second case (and symmetrically, the third) will follow from the second statement of the theorem, as the right-hand side of \eqref{E:tcol} is a corner.

\bigskip

The proof of the second statement requires some juggling with projections, so we pause here to sketch the idea.  First, if we specialize to the case where $\mathcal{M}$ and $\mathcal{N}$ are factors, \eqref{E:tcol} says that the image of a column is either a column or a row.  For non-factors and columns as described, the image is a central sum of a column and row, with $z'$ demarcating the two pieces.

The proof is effected by using the projection $q$ to divide $L^p(\mathcal{M})$ into four corners, each of which is an orthocomplement.  The original column is divided into two corners, $A$ and $B$, and we show that the ``checkerboard" array is preserved by $T$.  Visually,
$$L^p(\mathcal{M})q = \left( \begin{smallmatrix} A & 0 \\ B & 0 \end{smallmatrix} \right); \qquad L^p(\mathcal{M}) = \left( \begin{smallmatrix} A & B^\perp \\ B & A^\perp \end{smallmatrix} \right).$$
When the algebras are factors, there are only two (schematic) possibilities for $T$:
$$T: \left( \begin{smallmatrix} A & B^\perp \\ B & A^\perp \end{smallmatrix} \right) \mapsto \left( \begin{smallmatrix} T(A) & T(B)^\perp \\ T(B) & T(A)^\perp \end{smallmatrix} \right) \qquad \text{and} \qquad T: \left( \begin{smallmatrix} A & B^\perp \\ B & A^\perp \end{smallmatrix} \right) \mapsto \left( \begin{smallmatrix} T(A) & T(B) \\ T(B)^\perp & T(A)^\perp \end{smallmatrix} \right).$$
To show this we need to look hard at the pairs of projections defining $T(A)$ and $T(B)$.  We will see that either the left projections agree and the right projections are orthogonal with sum 1, or vice versa.  To make the bookkeeping a little more confusing, on non-factors the two possibilities can each happen on a central summand.

\bigskip

So we now assume the hypotheses of the second statement, and set $A = qL^p(\mathcal{M})q$, $B = (1-q) L^p(\mathcal{M})q$.  As argued in the fourth case above, $T(A)$ and $T(B)$ are corners, say $r_1 L^p(\mathcal{N})r_2$ and $s_1 L^p(\mathcal{N}) s_2$.  Since $A$ and $B$ neither contain nor are disjoint from any central summand, the same is true for $T(A)$ and $T(B)$.  It follows that $r_1,r_2, s_1, s_2$ are strictly between 0 and 1 on all central summands.

Substituting into \eqref{E:lpsip},
$$\xi \in A, \: \eta \in B \quad \Rightarrow \quad [\xi, \eta] = 0.$$
By Lemma \ref{T:preserve}, any pair of vectors from $T(A)$ and $T(B)$ also has semi-inner product zero.  Lemma \ref{T:corth} then tells us that the central supports $x_1$ of $r_1 s_1$ and $x_2$ of $r_2 s_2$ are orthogonal.

Notice that $T(B^\perp) = T(B)^\perp = (1-s_1)L^p(\mathcal{N})(1-s_2)$, and similarly for $T(A^\perp)$.  Now we apply the reasoning of the previous two paragraphs to the pair $A, B^\perp$, showing that the central supports $w_1$ of $r_1(1-s_1)$ and $w_2$ of $r_2 (1- s_2)$ are orthogonal.  But $w_1 \geq 1-x_1$, since
$$r_1 (1-s_1) (1-x_1)  = (r_1 - r_1 s_1) (1-x_1) = r_1 (1-x_1).$$
(The central support of the left-hand side is $\leq w_1$, of the right-hand side is $1-x_1$.)  Similarly $w_2 \geq 1-x_2$.  Since $x_1, x_2$ and $w_1, w_2$ are orthogonal pairs, we must have $w_1 = x_2$, $w_2 = x_1$, and $x_1 + x_2 = 1$. 

The preceding argument uses the pairs $(A, B)$ and $(A, B^\perp)$.  If we make the same argument for $(A, B)$ and $(A^\perp, B)$, then for $(A, B^\perp)$ and $(A^\perp, B^\perp)$, we may conclude that $x_1$ is the central support of each of $r_1 s_1, r_2(1-s_2), (1-r_2)s_2, (1-r_1) (1-s_1)$, while $x_2 = 1-x_1$ is the central support of each of $r_2 s_2, r_1(1-s_1), (1-r_1) s_1, (1-r_2)(1-s_2)$.  We write out two implications:
$$ r_2 s_2 x_1 = 0 = (1-r_2)(1-s_2) x_1 = (1 - r_2 - s_2 + r_2 s_2) x_1 \quad \Rightarrow \quad x_1 = (r_2 + s_2) x_1.$$
$$ (r_1 - r_1 s_1) x_1 = r_1 (1-s_1) x_1 = 0 = (1-r_1) s_1 x_1 = (s_1 - r_1 s_1) x_1 \Rightarrow r_1 x_1 = s_1 x_1.$$ 
Symmetrically $x_2 = (r_1 + s_1) x_2$ and $ r_2 x_2 = s_2 x_2 $.

Based on these last conclusions, we calculate
\begin{align*}
T(L^p(\mathcal{M})q) &= T(A) + T(B) \\
&= r_1 L^p(\mathcal{N})r_2 + s_1 L^p(\mathcal{N}) s_2 \\
&= (r_1 L^p(\mathcal{N})r_2 + s_1 L^p(\mathcal{N}) s_2) x_2 + (r_1 L^p(\mathcal{N})r_2 + s_1 L^p(\mathcal{M}) s_2) x_1 \\
&= L^p(\mathcal{N}) r_2 x_2 + r_1 L^p(\mathcal{N}) (1 - x_2), \\
\end{align*}
which verifies \eqref{E:tcol} by taking $q_1 = r_2$, $q_2 = r_1$, and $z' = x_2$.
\end{proof}

Note that the projections $q_1 z', q_2(1-z'), z'$ of Lemma \ref{T:allcorner}(2) are all uniquely determined by $q$.  Even more is true.

\begin{lem} \label{T:ind}
Assume that $\mathcal{M}$ has no abelian summand.  The central projection $z'$, defined in Lemma \ref{T:allcorner}(2), does not depend on the choice of $q$.
\end{lem} 
\begin{proof}
Of course, all choices are still assumed to be strictly between 0 and 1 on all central summands.  For projections other than $q$ we will use obvious variants of \eqref{E:tcol}.

First observe that $z'$ does not change if we replace $q$ by a smaller projection $\dot{q}$.  Just write
\begin{align*}
L^p(\mathcal{N}) q_1 z' + q_2 L^p(\mathcal{N}) (1-z') &= T(L^p(\mathcal{M})q) \\
&\supseteq T(L^p(\mathcal{M})\dot{q})\\
&= L^p(\mathcal{N}) \dot{q}_1 \dot{z}' + \dot{q}_2 L^p(\mathcal{N}) (1-\dot{z}').
\end{align*}
Since columns which contain no central summands never contain nonzero rows (and vice versa), we must have $z' = \dot{z}'$.

We also claim that $z'$ does not change if we replace $q$ by a projection $\ddot{q}$ with $q \wedge \ddot{q} = 0$.  In this case we get the disjointness of
$$T(L^p(\mathcal{M})q) = L^p(\mathcal{N}) q_1 z' + q_2 L^p(\mathcal{N}) (1-z') $$
and
$$T(L^p(\mathcal{M})\ddot{q}) = L^p(\mathcal{N}) \ddot{q}_1 \ddot{z}' + \ddot{q}_2 L^p(\mathcal{N}) (1-\ddot{z}').$$
A row and a column with overlapping central support always have nonzero intersection, so necessarily $z' = \ddot{z}'$.

Finally, given any other projection $r$, let $y$ be the central support of $q \wedge r$.  We may consider $T$ restricted to $L^p(\mathcal{M})y$; by Lemma \ref{T:corner} this is still an $L^p$ isometry.  The second paragraph shows that the (now restricted) projection $z'$ does not change if we go from $q$ to $q \wedge r$ to $r$.  Similarly, for $T$ restricted to $L^p(\mathcal{M})(1-y)$, the third paragraph allows us to pass from $q$ to $r$ without altering the restriction of $z'$.
\end{proof}

\begin{lem} \label{T:cccr}
Assume that $\mathcal{M}$ has no abelian summand, let $z'$ be as in Lemma \ref{T:ind}, and let $z$ be the corresponding central projection in $\mathcal{M}$.  Then on $L^p(\mathcal{M})z$, $T$ takes columns to columns, while on $L^p(\mathcal{M})(1-z)$, $T$ takes columns to rows.
\end{lem}

\begin{proof}
Let $L^p(\mathcal{M})r \subset L^p(\mathcal{M})z$ be a column containing no central summand, and let $z_0 \leq z$ be the central support of $r$.  Find a projection $\dot{r}$ with central support $(1-z_0)$ so that $L^p(\mathcal{M})(r + \dot{r})$ still contains no central summands.  Applying Lemma \ref{T:allcorner} for the projection $q = r + \dot{r}$,
\begin{align*}
T(L^p(\mathcal{M})r) &= T(L^p(\mathcal{M})(r + \dot{r})z_0) = [T(L^p(\mathcal{M})(r + \dot{r}))]z_0' \\
&= [L^p(\mathcal{N}) q_1 z' + q_2 L^p(\mathcal{N}) (1-z')] z_0' = L^p(\mathcal{N}) q_1 z_0'.
\end{align*}

An arbitrary column in $L^p(\mathcal{M})z$ is a central sum of a central summand and a column containing no central summands.  By Lemma \ref{T:corner} and the preceding paragraph, its image under $T$ is a central sum of columns, which is again a column.  The argument for $L^p(\mathcal{M})(1-z)$ is similar.
\end{proof}

Now we return to general $\mathcal{M}, \mathcal{N}$ and look to divide our problem into two pieces.  With $a \mathcal{M}$ the abelian summand of $\mathcal{M}$, we apply Lemma \ref{T:cccr} to the restriction $T:L^p(\mathcal{M})(1-a) \overset{\sim}{\to} L^p(\mathcal{N})(1-a')$.  This gives us a central projection $z \leq 1-a$ such that for any central projection $y$ with $z \leq y \leq z + a$, the restriction $T:L^p(\mathcal{M})y \overset{\sim}{\to} L^p(\mathcal{N})y'$ takes columns to columns, while $T:L^p(\mathcal{M})(1-y) \overset{\sim}{\to} L^p(\mathcal{N})(1-y')$ takes columns to rows.  (On $L^p(\mathcal{M})a$ and $L^p(\mathcal{N})a'$, there is no difference between columns, rows, and corners, as all are central summands.)  For now we focus on one piece, renaming $L^p(\mathcal{M})y$ as $L^p(\mathcal{M})$, $L^p(\mathcal{N})y'$ as $L^p(\mathcal{N})$, and the restriction of $T$ as $T$.  We have that
\begin{equation} \label{E:rtsupp}
T(L^p(\mathcal{M}) q) = L^p(\mathcal{N})\pi_r(q)
\end{equation}
for a well-defined increasing map $\pi_r$ between projection lattices.

It follows from Lemma \ref{T:allcorner}(2) that when $q$ is strictly between 0 and 1 on every central summand which contains no abelian summand, $\pi_r(q)$ is as well.  So if we apply Lemmas \ref{T:allcorner}, \ref{T:ind}, \ref{T:cccr} to $T^{-1}$, we see that $T^{-1}$ also takes columns to columns, and both $T$ and $T^{-1}$ take rows to rows.  More importantly, $\pi_r$ is bijective.  Now equation \eqref{E:rtsupp} implies that $s_r(T(\xi)) \leq \pi_r(s_r(\xi))$ for any $\xi \in L^p(\mathcal{M})$.  Since we can make the same argument for $T^{-1}$, we must actually have that
\begin{equation} \label{E:rsupp}
s_r(T(\xi)) = \pi_r(s_r(\xi)), \qquad \xi \in L^p(\mathcal{M}).
\end{equation}  

We claim that $\pi_r$ preserves orthogonality of projections.  Indeed, if $e \perp f$ in $\mathcal{P}(\mathcal{M})$, then any $\xi \in L^p(\mathcal{M})e$ and $\eta \in L^p(\mathcal{M})f$ have semi-inner product zero.  Combining Lemmas \ref{T:preserve}, \ref{T:corth}, and equation \eqref{E:rtsupp} gives $\pi_r(e) \perp \pi_r(f)$. 

Let $\xi \in L^p(\mathcal{M})$ and $p \in \mathcal{P}(\mathcal{M})$.  Using \eqref{E:rsupp} and properties of $\pi_r$,
\begin{equation} \label{E:rmodp}
T(\xi p) = T(\xi p) \pi_r(p) = T(\xi p) \pi_r(p) + T(\xi (1-p)) \pi_r(p) = T(\xi) \pi_r(p).
\end{equation}
Now we extend $\pi_r$ in a standard way: first by linearity to real linear combinations of orthogonal projections, then by continuity to self-adjoint elements, then by the equation
$$\pi_r(x + iy) = \pi_r(x) + i\pi_r(y), \qquad x,y \in \mathcal{M}_{sa},$$
to all of $\mathcal{M}$.  To see that $\pi_r$ is linear, note that by construction we have
\begin{equation} \label{E:rmod}
T(\xi x) = T(\xi) \pi_r(x), \qquad \forall \xi \in L^p(\mathcal{M}), \forall x \in \mathcal{M}.
\end{equation}
So for any $x,y \in \mathcal{M}$, $\xi \in L^p(\mathcal{M})$,
$$T(\xi) (\pi_r(x) + \pi_r(y)) = T(\xi x) + T(\xi y) = T(\xi(x+y)) = T(\xi) \pi_r(x+y),$$
which implies $\pi_r(x) + \pi_r(y) = \pi_r(x+y)$.  By construction $\pi_r$ is *-preserving and bijective.  Finally, take $x,y \in \mathcal{M}$, $\xi \in L^p(\mathcal{M})$, and calculate 
\begin{equation} \label{E:mult}
T(\xi) \pi_r(xy) = T(\xi xy) = T(\xi x)\pi_r(y) = T(\xi)\pi_r(x)\pi_r(y).
\end{equation}
Apparently $\pi_r: \mathcal{M} \to \mathcal{N}$ is also multiplicative.

Being a surjective *-isomorphism, $\pi_r$ induces a surjective isometry from $L^p(\mathcal{M})$ to $L^p(\mathcal{N})$ as discussed earlier in this section. We will denote this map by $\rho$: key properties are
$$\rho(x \xi y) = \pi_r(x) \rho(\xi) \pi_r(y), \qquad \rho(\varphi^{1/p}) = (\varphi \circ \pi_r^{-1})^{1/p},$$
for $x,y \in \mathcal{M},$ $\xi \in L^p(\mathcal{M}),$ $\varphi \in \mathcal{M}_*^+$.

Now we consider the surjective isometry $T \circ \rho^{-1}: L^p(\mathcal{N}) \to L^p(\mathcal{N})$.  This is actually a right module map:
$$T \circ \rho^{-1}(\xi x) = T(\rho^{-1}(\xi) \pi_r^{-1}(x)) = T \circ \rho^{-1}(\xi) x, \qquad x \in \mathcal{N}, \: \xi \in L^p(\mathcal{N}).$$
It is known that the left and right module actions of $\mathcal{N}$ on $L^p(\mathcal{N})$ are commutants of each other.  (This was first shown in \cite[Proposition 35]{T1}, or see \cite[Theorem 1.5]{JS} for a stronger result.)  Thus $T \circ \rho^{-1}$ is given by left multiplication by an element of $\mathcal{N}$, and by \cite[Lemma 1.1]{JS} the element has norm equal to $\|T\| = 1$.  The same is true for $[T \circ \rho^{-1}]^{-1}$, so the element is unitary - call it $u$.  Then for all $\varphi \in \mathcal{M}_*^+$,
\begin{equation} \label{E:concl}
T(\varphi^{1/p}) = T \circ \rho^{-1} \circ \rho (\varphi^{1/p}) = u \rho(\varphi^{1/p}) = u (\varphi \circ \pi_r^{-1})^{1/p}, 
\end{equation}
which was to be shown.

What about the case where $T$ takes columns to rows and vice versa?  Equation \eqref{E:rsupp} becomes $s_r(T(\xi)) = \pi_r(s_\ell(\xi))$, \eqref{E:rmod} becomes $T(x \xi) =  T(\xi) \pi_r(x)$, and a calculation parallel to \eqref{E:mult} shows that $\pi_r$ is \textit{antimultiplicative}.  Associating the $L^p$ isometry $\rho$ to $\pi_r$ as before, we still have that $T \circ \rho^{-1}$ is a right module map, and the conclusion \eqref{E:concl} follows.  So in the general case with both summands present, we may take the sum of the two partial isometries as the unitary $w$, and the sum of the *-homomorphism and *-antihomomorphism as the surjective Jordan *-isomorphism $J$.  The proof of Theorem \ref{T:typical} is complete.

\section{Remarks on the proof} \label{S:remarks}

\S 1. We chose to work with columns in $L^p(\mathcal{N})$ because of the desired polar decomposition.  In the multiplicative case handled first, one can also find a surjective *-isomorphism $\pi_\ell: \mathcal{M} \to \mathcal{N}$ such that equation \eqref{E:rmod} becomes $T(y \xi x) = \pi_\ell(y) T(\xi) \pi_r(x)$.  Moreover we have $\pi_\ell(y) = u \pi_r(y) u^*.$  Obvious variants hold for the antimultiplicative and general cases.

\bigskip

\S 2. It is possible to obtain the main results of this paper without using semi-inner products.  There is an alternate route to Lemma \ref{T:allcorner} which is essentially simpler, but unfortunately it does not apply to algebras with finite type I summands.  So in order to build a complete proof in this way, one must also isolate the finite type I summands by methods similar to Lemma \ref{T:abelian}, and apply there a known result (like  \cite[Theorem 2]{Y}).  We found it preferable to give a unified proof, with no dependence on type or previous isometry results.

However, the alternate proof has the significant advantage of applying equally well to $0 < p \leq 1$.  (For $0 < p < 1$, $L^p(\mathcal{M})$ is a $p$-Banach space.)  Since this may be of interest to some readers, the argument is featured in \cite{S3}, where the case $p=1$ is carried out explicitly and used to give a new proof of the noncommutative Banach-Stone theorem.  (By this we mean the nonunital C*-algebra version of Theorem \ref{T:kadison}, which was first stated by Paterson and Sinclair \cite{PS} in 1972.)  Therefore Theorems \ref{T:main} and \ref{T:typical} are also true for $0 < p \leq 1$.  The only other amendment to their proofs is that equations \eqref{E:alphap} and \eqref{E:bimod} must be justified directly, as interpolation cannot be used.

\bigskip

\S 3. Equation \eqref{E:rsupp} is already enough to settle the implication (2) $\to$ (1) in Theorem \ref{T:main}.  A bijective map between projection lattices which preserves orthogonality is called an \textbf{orthoisomorphism}; Dye \cite{D} showed that such a map is the restriction of a surjective Jordan *-isomorphism off of the type $\text{I}_2$ summand.  Since $\pi_r^{-1}$ is also an orthoisomorphism, it follows that $T$ maps the $\text{I}_2$ summands to each other.  But $T$ induces an isomorphism of centers, so the $\text{I}_2$ summands have isomorphic centers and are therefore also *-isomorphic.

\bigskip

\S 4. One can use Lemma \ref{T:allcorner}(1) to define the following map:
\begin{equation} \label{E:notorth}
(q_1, q_2) \mapsto (S_\ell(q_1, q_2), S_r(q_1, q_2)), \qquad q_1,q_2 \in \mathcal{P}(\mathcal{M}),
\end{equation}
where $S_\ell(q_1, q_2)$ and $S_r(q_1, q_2)$ are the unique projections in $\mathcal{N}$ with identical central support satisfying
$$T(q_1 L^p(\mathcal{M}) q_2) = S_\ell(q_1,q_2)L^p(\mathcal{N})S_r(q_1, q_2).$$
Because $T$ preserves orthogonality, the map \eqref{E:notorth} is almost an orthoisomorphism from $\mathcal{M} \oplus \mathcal{M}$ to $\mathcal{N} \oplus \mathcal{N}$.  The deficit has to do with central support; if one requires that the two inputs have identical central support, \eqref{E:notorth} ``densely defines" an orthoisomorphism.  In fact it is possible to show the strong continuity of this map and in this way construct an actual orthoisomorphism, at least when $\mathcal{M}$ has no type I summand.

Edwards and R\"{u}ttimann \cite{ER} specifically studied $CP(\mathcal{M})$, the complete *-lattice of pairs of projections with equal central support.  Just as we have suggested this as a tool for studying $L^p$ corners, they use an equivalence with the set of $L^\infty$ corners $q_1 \mathcal{M} q_2$.  This, in turn, is naturally equivalent to the lattice of structural projections and the lattice of weak*-closed inner ideals, both defined in terms of the Jordan triple structure of $\mathcal{M}$.  Their paper actually formalizes, in the language of lattice theory and Jordan triple systems, some of our manipulations of corners.

\bigskip

\textbf{Acknowledgments.} The author would like to thank Marius Junge for helpful input, Keiichi Watanabe for sharing his preprints, and the anonymous referee for several useful suggestions.

\end{document}